\documentclass[11pt,reqno]{amsart}

\usepackage{amsthm}
\usepackage{amssymb}
\usepackage{latexsym}
\usepackage{multicol}
\usepackage{verbatim,enumerate}
\usepackage{accents}
\usepackage[usenames]{color}
\usepackage{youngtab}
\usepackage[all, poly]{xy}  
\usepackage[utf8]{inputenc}
\usepackage{hyperref}
\usepackage{amsmath, amscd}
\usepackage{soul}
\usepackage{tikz} 
\usetikzlibrary{matrix,arrows,decorations.pathmorphing}

\theoremstyle{plain}
\newtheorem*{prop}{Proposition}
\newtheorem{thm}{Theorem}
\newtheorem*{lem}{Lemma}
\newtheorem*{cor}{Corollary}

\theoremstyle{definition}
\newtheorem*{example}{Example}
\newtheorem*{defn}{Definition}

\newtheorem*{rem}{Remark}

\theoremstyle{remark}

\newcommand{\lie}[1]{\mathfrak{#1}}

\newcommand\bn{\mathbb N}
\newcommand\bz{\mathbb Z}

\newcommand{\bs}{\mathbf{S}}

\def\a{\alpha}

\advance\textwidth by 1.2in  \advance\oddsidemargin by -.6in \advance\evensidemargin by -.6in

\newenvironment{pf}{\proof}{\endproof}
\newcounter{cnt}
 \makeatletter
\def\mydggeometry{\makeatletter\dg@YGRID=1\dg@XGRID=20\unitlength=0.003pt\makeatother}
\makeatother \theoremstyle{remark}

\numberwithin{equation}{section}
\makeatletter
\def\section{\def\@secnumfont{\mdseries}\@startsection{section}{1}%
  \z@{.7\linespacing\@plus\linespacing}{.5\linespacing}%
  {\normalfont\scshape\centering}}
\def\subsection{\def\@secnumfont{\bfseries}\@startsection{subsection}{2}%
  {\parindent}{.5\linespacing\@plus.7\linespacing}{-.5em}%
  {\normalfont\bfseries}}
\makeatother

\begin{document}

\title[Simplified presentations and embeddings of Demazure modules]{Simplified presentations and embeddings of Demazure modules}
\author{Deniz Kus}
\address{University of Bochum, Faculty of Mathematics, Universit{\"a}tsstr. 150, 44801 Bochum, 
Germany}
\email{deniz.kus@rub.de}
\author{R. Venkatesh}
\address{Department of Mathematics, Indian Institute of Science, Bangalore 560012}
\email{rvenkat@iisc.ac.in}
\thanks{D.K. was partially funded by the Deutsche Forschungsgemeinschaft (DFG, German Research Foundation)-- grant 446246717. R.V. was partially funded by the grants DST/INSPIRE/04/2016/000848, MTR/2017/000347, and an Infosys Young Investigator Award.}
\begin{abstract}
For an untwisted affine Lie algebra we prove an embedding of any higher level Demazure module into a tensor product of lower level Demazure modules (e.g. level one in type A) which becomes in the limit (for anti-dominant weights) the well-known embedding of finite-dimensional irreducible modules of the underlying simple Lie algebra into the tensor product of fundamental modules. To achieve this goal, we first simplify the presentation of these modules extending the results of \cite{CV13} in the $\lie g$-stable case. As an application, we propose a crystal theoretic way to find classical decompositions with respect to a maximal semi-simple Lie subalgebra by identifying the Demazure crystal as a connected component in the corresponding tensor product of crystals.
\end{abstract}

\maketitle

\section{Introduction}
Let $\lie g$ be a simple finite-dimensional complex Lie algebra with corresponding untwisted affine Lie algebra $\widehat{\lie g}$. Affine Demazure modules are a family of finite-dimensional representations that are characterized by an element in the affine Weyl group and a dominant integral weight. They have been the subject of intense study with lots of applications and conjectures (see for example \cite{CV13, FoL07,I03,Kumar, KV14, Na11,San00} and the references therein). Although there are character formulas in terms of Demazure operators \cite{Kumar}, it seems very complicated to give closed dimension formulas let alone closed character formulas which are comparable to the Kac-Weyl character formula. 

The first connection to Macdonald polynomials was observed in \cite{San00} and later extended to the simply-laced case in \cite{I03}, where the level one Demazure characters are identified with certain specializations of Macdonald polynomials. For non simply-laced types a similar approach connects the specialized Macdonald polynomials with (generalized) Weyl modules \cite{CI15,FM17a}. In the higher level case character or dimension formulas for Demazure modules are generically not known and this problem has been partially solved for $\mathfrak{sl}_{n+1}$-stable Demazure modules of level two in \cite{BK20b,Biswal-Chari}. We emphasize at this point that the most critical ones are the prime Demazure modules which do not admit a non-trivial splitting into a tensor product. 

The article \cite{BK20b} uses the dual functional realization of the loop algebra and the fact that level two Demazure modules for $\mathfrak{sl}_{n+1}$ appear in a completely different context. They are the graded limits of an important family of modules for the quantum affine algebra and the classical decomposition is obtained via lattice points in convex polytopes. The article \cite{Biswal-Chari} uses a different approach which was first introduced and developed in \cite{Jos06}, namely the theory of Demazure flags.  This approach provides a further deep and unexpected connection to combinatorics and number theory; see for example \cite{BK19a} for the connection to the combinatorics of Dyck path or \cite{BCK18} and \cite{BCSV15} for the connection to Mock-Theta functions and hypergeometric series. 

Another approach to the problem of finding dimension formulas, character formulas or classical decompositions is given by crystal theory. It is known that the Demazure module $V_w(\Lambda)$ admits a crystal base \cite{Ka93a} which is the full subgraph of the crystal base of the irreducible integrable highest weight module whose vertices consist precisely of those elements that are reachable by raising operators from the unique element of weight $w\Lambda$. There is a connection to the tensor product of Kirillov-Reshetikhin crystals which has been worked out by various authors in the last years (see for example \cite{LS19a,Naoi,ST12} and the references therein). The main result of \cite{LS19a}, which seems to be the most general version so far, states that the $\mathfrak{g}$-stable level $k$ Demazure module for non-exceptional types appear inside a tensor product of Kirillov-Reshetikhin crystals by removing all edges that are not level $k$ Demazure edges (see \cite[Section 2]{LS19a} for a precise definition) and picking the connected component which contains the unique element of extremal weight $w\Lambda$. Finding very explicit combinatorial models for these crystals is still a very important and open problem. 

Most of the connections described above benefit enormously from the presentation of Demazure modules proved by Mathieu in \cite{M88}. However the relations are sometimes very hard to check. After simplifying the defining relations of $\lie g$-stable Demazure modules in \cite{CV13} many new isomorphisms and maps were proved afterwards. By way of example, in the same article, the authors proved that the Q-system introduced by Kirillov and Reshetikhin \cite{KR87a} extends to a canonical short exact sequence of fusion products. In the study of global Demazure modules and arc schemes \cite{DF19a} the authors used the simplified presentation to show that rectangular $\lie g$-stable Demazure modules can be obtained as localizations of global Demazure modules at zero. 
The same simplification was used in \cite{KL14} to show that the fusion product of finite-dimensional irreducible $\lie g$-modules of the same highest weight is isomorphic to the truncated Weyl module or in \cite{BCMo15} to identify level two $\lie g$-stable Demazure modules with graded limits of important families of modules for the quantum affine algebra. 

The crystal theoretic approach to Demazure modules and the realization inside a tensor product obtained in \cite{LS19a} is one motivation of the article to embed any higher level Demazure crystal (no restrictions on the type nor on the extremal weight) into a tensor product of Demazure crystals which are in some way easier to understand. We attack this problem algebraically which we will explain now in more details. 

For an element $w$ in the affine Weyl group and a dominant integral weight $\Lambda$ let $V_w(\Lambda)=\mathbf{D}_{\mu}^k$ the affine Demazure module where $w\Lambda=\mu+k\Lambda_0$ (modulo the null root). For such an element $\mu\in P$ we define the set 
$$\mathcal{P}(\mu,k)=\{(\mu_1,\dots,\mu_k)\in P^{\times k}: \mu_1+\dots+\mu_k=\mu\}$$
and introduce the notion of $r$-admissible elements in Section~\ref{section4}. We show that for an $r$-admissible element $(\mu_1,\dots,\mu_k)$ (sometimes the unique element described in \cite{F14aa} is the correct choice) we have an embedding of graded modules
$$\mathbf{D}_{\mu}^{rk}\hookrightarrow \mathbf{D}_{\mu_1}^{r}\otimes \cdots \otimes \mathbf{D}_{\mu_k}^{r}.$$
It is clear that this embedding becomes in the limit $k\rightarrow \infty$  for anti-dominant weights the well-known embedding of finite-dimensional irreducible modules of the underlying simple Lie algebra into the tensor product of fundamental modules. To achieve this embedding we first need to simplify the presentation of Demazure modules for
any weight $\mu$ (see Section~\ref{section41}) extending the $\lie g$-stable case from \cite{CV13}. We believe that this fact is of independent interest and finds several more applications in the future.

Finally we propose, similar to the approach taken in \cite{LS19a}, a way to come up with classical decompositions. For any weight $\mu\in P$ we define a subset $R^{-}(\mu)$ (see Section~\ref{section2}) and if this set is non-empty we can find a maximal semi-simple subalgebra $\lie g_0\subseteq\lie g$ corresponding to the nodes in $R^-(\mu)$; the $\lie g$-stable case corresponds to $R^+=R^-(\mu)$. The $\lie g_0$ decompostion of the Demazure module ``seems'' (we use this wording because it is not made precise in this article) to be encoded in the connected component containing the unique element of weight $w\Lambda$ inside the tensor product (according to the splitting above), after removing the arrows not corresponding to the nodes in $R^-(\mu)$. The problem of describing the corresponding crystal inside the tensor product in a very explicit way will be discussed elsewhere.\\\

\textit{Organization of the paper:} In Section~\ref{section2} we establish the basic notation and elementary results needed in the rest of the paper. In Section~\ref{section3}, we define three classes of modules $M_{\mu,\mathbf{p}},M'_{\mu,\mathbf{p}}$ and $M''_{\mu,\mathbf{p}}$, $\mu$ an arbitrary integral weight and $\mathbf{p}$ a sequence of functions indexed by the positive roots, and prove surjective maps among them. This allows us to simplify the presentation of all Demazure modules in Section~\ref{section4} which was obtained earlier in \cite{CV13} in the $\lie g$-stable case. Furthermore, we prove an embedding of Demazure modules and propose a way to calculate the classical decompositions with respect to a maximal semi-simple Lie subalgebra. \\\

\textit{Acknowledgement:} The first author thanks Daisuke Sagaki for many helpful discussions.
\section{Preliminaries}\label{section2}

\subsection{} We denote by $\mathbb{Z}$, $\mathbb{Z}_+$, and $\mathbb{N}$ the set of integers, the set of non-negative integers and the set of positive integers respectively. The base field will be the complex numbers $\mathbb{C}$ throughout. 
For a given Lie algebra $\lie a$, we denote by $\mathbf{U}(\lie a)$ the universal enveloping algebra of $\lie a$ and by $\lie a[t]=\lie a\otimes \mathbb{C}[t]$ its current algebra with Lie bracket
$$[x\otimes t^r, y\otimes t^s]=[x,y]\otimes t^{r+s},\ \ x,y\in\lie a,\ \ r,s\in \mathbb{Z}_+.$$
\subsection{} Let $\mathfrak{g}$ be a simple finite-dimensional complex Lie algebra with Cartan subalgebra $\mathfrak{h}$ and root system $R$ spanned by the simple roots $\{\alpha_1,\dots,\alpha_n\}$ with index set $I=\{1,\dots,n\}$. Let $\mathfrak{b}$ be a Borel subalgebra and $R^{\pm}$ be the corresponding set of positive and negative roots respectively and denote by $\theta\in R^+$ the highest root of $\lie g$. Let $(\cdot,\cdot)$ be the non-degenerate bilinear form on $\lie h^{*}$ with $(\theta,\theta)=2$ induced by the restriction of the (suitably normalized) Killing form of $\lie g$ to $\lie h$ and set $d_{\alpha}=2/(\alpha,\alpha)$ for $\alpha\in R$. We denote by $P$ the set of integral weights with basis $\{\varpi_1,\dots,\varpi_n\}$ and let $P^+$ be the subset of dominant integral weights. We further fix a Chevalley basis $\{x_{\alpha}^{\pm}, h_{\alpha}: \alpha\in R^+\}$ of $\mathfrak{g}$ and set 
$$R^{+}(\mu)=\{\alpha\in R^+: \mu(h_{\alpha})\ge 0\},\ \ R^-(\mu)=\{\alpha\in R^+: \mu(h_{\alpha})\leq 0\},\ \ \ \mu\in P.$$
We denote by $\widehat{\lie g}$ the corresponding affine Kac-Moody algebra with underlying simple Lie algebra $\mathfrak{g}$. The canonical central element is denoted by $c$ and the derivation by $d$ and let $\widehat{\lie b}$ be the standard Borel subalgebra. The induced set of positive and negative roots respectively is denoted by $\widehat{R}^{\pm}$ and let $\widehat{P}$ and $\widehat{P}^+$ respectively the set of integral and dominant integral weights respectively. Let $W$ be the Weyl group of $R$ and recall that the affine Weyl group $\widehat{W}$ is given by the semi-direct product
$$\widehat{W}=W\ltimes t_{M},\ \ t_{M}=\{t_\mu: \mu\in M\}$$
 where $M$ is the lattice  generated by the elements $w(\theta)$ for all $w\in W$ (see \cite[Proposition 6.5]{K90}). 

 \subsection{}
Given $\Lambda\in \widehat{P}^+$, let $V(\Lambda)$ be the irreducible, integrable $\widehat{\lie g}$-module with highest weight $\Lambda$ and highest weight vector $v_\Lambda$. We have a weight space decomposition $$V(\Lambda)=\bigoplus_{\gamma} V(\Lambda)_{\Lambda-\gamma},\ \ \ \mathrm{dim} V(\Lambda)_{w\Lambda}=1,\ \ \forall w\in \widehat{W}.$$ An important class of representations for the Borel subalgebra is given by the so-called \textit{Demazure modules}. To each pair
$(\Lambda, w)\in \widehat{P}^+\times  \widehat{W}$ we define $V_w(\Lambda)=\mathbf{U}(\widehat{\lie b})V(\Lambda)_{w\Lambda}$ and note that
$$w\Lambda=w'\Lambda \Rightarrow V_w(\Lambda)= V_{w'}(\Lambda).$$ 
We denote these modules alternatively by $\mathbf{D}^k_{\mu}[i]$ whenever 
\begin{equation}\label{zzz}w\Lambda=\mu+k \Lambda_0+i \delta\end{equation} 
where $\Lambda_0\in \widehat{\mathfrak{h}}^*$ is the $0$-th fundamental weight determined by
$$\Lambda_0(c)=1,\ \Lambda_0(h_{\alpha_i})=0=\Lambda_0(d),\ \ i=1.\dots,n.$$
Moreover, for each triple $(\mu,k,i)$ there exists a Weyl group element $w$ and a dominant affine weight $\Lambda$ such that \eqref{zzz} holds. This follows immediately from the fact that $P+\mathbb{N}\Lambda_0+\mathbb{Z}\delta$ lies in the Tits cone \cite[Proposition 5.8]{K90}.
\begin{rem} The modules $\mathbf{D}^k_{\mu}[i]$ are $\lie g$-stable if and only if $\mu$ is an anti-dominant weight, and in this case they are often denoted by $D(k,w_0\mu)[i]$ in the literature where $w_0$ is the longest element in the Weyl group $W$ (see for example \cite{BCMo15, CV13, Na11}). However $D(k,w_0\mu)[i]$ or $D_{k,w_0\mu}[i]$ is also used sometimes to denote the module $\mathbf{D}_k^{k\mu}[i]$ of lowest weight $k\mu$ (see for example \cite{DF19a}) which leads to lot of confusion. We decided to use the notation introduced above to avoid any confusion with the literature and to emphasize that we allow arbitrary weights.
\end{rem}

A presentation for Demazure modules has been proved in \cite{M88}, which we record now in the affine setting. In the finite--dimensional setting this result was obtained earlier in \cite{Jo85a}.
\begin{thm}\label{mathieu}
The module $\mathbf{D}^k_{\mu}[i]$ is as a $\mathbf{U}\big(\widehat{\mathfrak{b}}\big)$--module isomorphic to the cyclic module generated by a non-zero vector $v$ with the following defining relations: For  $h\in \mathfrak{h}$ and
$\alpha\in R^+$, we have
\begin{itemize}
\item $(h\otimes t^s)v=\delta_{s,0}\mu(h)v$, $dv=i v$, and $cv=kv$,\vspace{0,1cm}
\item $(x_\alpha^{\pm}\otimes t^{s^{\pm}})^{p^{\pm}_{\alpha}+1}v=0$,\ \ $p^{\pm}_{\alpha}=\max\{0,\mp\mu(h_{\alpha})-d_{\alpha}s^{\pm}k\}$,\ \ $s^+\geq 0$,\ \ $s^->0.$
\end{itemize}
\end{thm}
The above presentation of these modules has been vastly simplified in \cite{CV13} (see \cite{KV14} for the twisted case) when $\mu$ is anti-dominant. This simplification turned out to be quite useful in proving the existence of maps from Demazure modules to other classes of modules. By way of example see \cite[Section 5]{CV13} where the authors proved the fact that the Q-system introduced by Kirillov and Reshetikhin \cite{KR87a} extends to a canonical short exact sequence of fusion products or in the study of global Demazure modules and arc schemes \cite{DF19a} where the authors showed that rectangular Demazure modules can be obtained as localizations of global Demazure modules at zero.
For more applications of the simplified presentation we refer the reader to \cite{FMM19a,KL14,Ra14a,Ve15a}.

\section{The main modules and existence of maps}\label{section3}
 One motivation of the paper is to simplify the defining relations of a class of representations for the Iwahori subalgebras $\mathfrak{I}[t]:=\mathfrak{b}\oplus (\mathfrak{g}\otimes t  \mathbb{C}[t])$ in untwisted affine Kac-Moody algebras. This class will include in particular all Demazure modules and we focus on the simplification of them in the next section. As an application we will prove an embedding of any higher level Demazure module into a tensor product of lower level Demazure modules (see Section~\ref{section4}). This allows us to obtain lower bounds for the dimensions of certain irreducible modules for the quantum affine algebra in simply-laced type.
 \subsection{}\label{section31} We first introduce this class of representations. Let $\mathbf{p}=(p_{\alpha}^{\pm})$ be a sequence of functions indexed by the positive roots
$$p^+_{\alpha}: \mathbb{Z}_+\rightarrow \mathbb{Z}_+,\ \ p^-_{\alpha}: \mathbb{N}\rightarrow \mathbb{Z}_+,\ \ \alpha\in R^+$$
such that there exists a minimal $s^{\pm}_{\alpha}$ satisfying $p^{\pm}_{\alpha}(a)=0$ for all $a\geq s^{\pm}_{\alpha}$. Given $\mu\in P$ we sometimes extend without further comment the function $p_{\alpha}^-$ to the non-negative integers by setting $p_{\alpha}^-(0)=\max\{0,\mu(h_{\alpha})\}$.
\begin{defn}\label{maindef} 
For an integral weight $\mu\in P$ we denote by $M_{\mu,\mathbf{p}}$ the cyclic module for the Iwahori subalgebra $\mathfrak{I}[t]$ generated by a non-zero vector $v$ with the following defining relations. For $\alpha\in R^{\mp}(\mu)$ and $h\in\mathfrak{h}$  we have:\vspace{0,1cm}
\begin{enumerate}
\item $\big(h\otimes t^s\big ) v=\delta_{s,0} \mu(h)v$,\ \ $s \geq 0$,\ \  $\big(x_\alpha^{\mp}\otimes t^{s^{\mp}}\big)v=0,\ \ s^+\geq 0,\ \ s^->0$ \vspace{0,1cm}

\item $\big(x_\alpha^{\pm}\otimes t^{\epsilon^{\pm}}\big)^{p^{\pm}_\alpha(\epsilon^{\pm})+1}v$, where $\epsilon^+=0$ and $\epsilon^-=1$\vspace{0,1cm}

\item\label{wegl} For all $1\leq i\leq s_{\alpha}^{\pm}$ and tuples $(a_i,..., a_{s_{\alpha}^{\pm}})$ satisfying $\displaystyle\sum_{k=0}^{s_{\alpha}^{\pm}-i}(k+1)a_{i+k}\geq p^{\pm}_\alpha(i)+1$ we have 
$$\big(x_\alpha^{\pm}\otimes t^{s_{\alpha}^{\pm}}\big)^{a_{s_{\alpha}^{\pm}}}\big(x_\alpha^{\pm}\otimes t^{s_{\alpha}^{\pm}-1}\big)^{a_{s_{\alpha}^{\pm}-1}}\cdots \big(x_\alpha^{\pm}\otimes t^{i}\big)^{a_i}v=0.$$
\end{enumerate}
\end{defn}
These modules are needed later to simplify the relations of the class of modules $M''_{\mu,\boldsymbol{p}}$ which we define now.
\subsection{}Given $\mu\in P$ and $\mathbf{p}$ as above, consider the tuple \begin{equation}\label{part1}_{\pm}\boldsymbol{\xi}^{\alpha}=(_{\pm}\xi_{1}^{\alpha},_{\pm}\xi_{2}^{\alpha},\dots),\ \ _{\pm}\xi_{i}^{\alpha}=p_{\alpha}^{\pm}(i-1)-p_{\alpha}^{\pm}(i), \ \alpha\in R^{\mp}(\mu)\end{equation}
\begin{defn}\label{maindef2} Define $M'_{\mu,\mathbf{p}}$ to be the cyclic $\mathfrak{I}[t]$-module with defining relations given by Definition~\ref{maindef}(1)-(2) and imposing the relations in Definition~\ref{maindef}(3) (for fixed $i$) only if 
$$_{\pm}\xi_{i+1}^{\alpha}< \ _{\pm}\xi_{i}^{\alpha}\geq \text{ $a_i+\cdots+a_{s_{\alpha}^{\pm}}$}.$$
Moreover, denote by $M''_{\mu,\mathbf{p}}$ to be the cyclic $\mathfrak{I}[t]$-module with defining relations given by Definition~\ref{maindef}(1)-(2) and 
$$(x_{\alpha}^{\pm}\otimes t^i)^{p_{\alpha}^{\pm}(i)+1},\ \ \alpha\in R^{\mp}(\mu),\ \ 1\leq i\leq s_{\alpha}^{\pm}$$
\end{defn}
For example, we always include the relation $(x_{\alpha}^{\pm}\otimes t^{s_{\alpha}^{\pm}})v=0$ in $M'_{\mu,\mathbf{p}}$ and $M''_{\mu,\mathbf{p}}$ since $$0=\ _{\pm}\xi_{s^{\pm}_{\alpha}+1}^{\alpha}< \ _{\pm}\xi_{s_{\alpha}^{\pm}}^{\alpha}\geq 1$$ by the minimality of $s_{\alpha}^{\pm}$. 
The importance of $M''_{\mu,\boldsymbol{p}}$ will become clear by the following example. They appear as many well-studied representations hence giving a strong connection to combinatorics, graded limits of representations for quantum affine algebras \cite{BCMo15,CP01} and the theory of Macdonald polynomials \cite{CI15,FM17a,I03}.
\begin{example}\label{mainex} We discuss three examples.
\begin{enumerate}
\item  Let $\mu\in P$ be an anti-dominant weight. If we set $p_{\alpha}^-(s)=0$ for all $s\in\mathbb{N}$ and $p_{\alpha}^+(s)=\max\{0,-\mu(h_{\alpha})-s\}$ for all $s\in\mathbb{Z}_+$ we recover the local Weyl module, namely we have an isomorphism
 $W_{\mathrm{loc}}(w_0\mu)\cong M_{\mu,\mathbf{p}}''$ as $\lie g[t]$-modules; we can extend the $\mathfrak{I}[t]$-structure on $M_{\mu,\mathbf{p}}''$ by requiring that $(\mathfrak{n}^-\otimes 1)$ acts as zero on $v$.
\item  Let $\mu\in P$ and $k\in\mathbb{N}$. We set 
$$p_{\alpha}^{\pm}(s)=\mathrm{max}\{0, \mp\mu(h_{\alpha})-d_\alpha s k\},\ \ s\in\mathbb{N},\ \ p_{\alpha}^{+}(0)=\mathrm{max}\{0, -\mu(h_{\alpha})\}$$
Then, by Theorem~\ref{mathieu} we have an isomorphism $M''_{\lambda,\mathbf{p}}\cong \mathbf{D}_{\mu}^{k}$, where we drop the dependence on $i$ since $\mathbf{D}_{\mu}^{k}[i]\cong \mathbf{D}_{\mu}^{k}[i']$ as $\mathfrak{I}[t]$-modules.
\item Choosing the data
$$p_{\alpha}^{+}(0)=0,\ p_{\alpha}^{-}(1)=\mu(h_{\alpha}),\ \ \alpha\in R^{+}(\mu)$$
and 
$$p_{\alpha}^{+}(0)=-\mu(h_{\alpha}),\ p_{\alpha}^{-}(1)=0,\ \ \alpha\in R^{-}(\mu)$$
we get an isomorphism to the generalized Weyl modules studied in \cite{FM17a}; in fact we have to choose the remaining integers as well, but we can simply choose the remaining $p_{\alpha}^{\pm}(s)$ so that they become redundant.
\end{enumerate}
\end{example}
The following simple observation is needed later.
\begin{lem}\label{n234}
Let $\mathbf{p}$ be the sequence of functions from Example~\ref{mainex}(2). Then, we have 
\begin{equation}\label{tosim}2 p_{\alpha}^{\pm}(i)\leq p_{\alpha}^{\pm}(i+1)+p_{\alpha}^{\pm}(i-1)\end{equation}
and equality holds if one of the following two conditions hold\begin{itemize}
\item $1\leq i\leq s_{\alpha}^{\pm}-2$ 
\vspace{0,15cm}

\item $i= s_{\alpha}^{\pm}-1$ and $\pm\mu( h_{\alpha})=d_{\alpha}s_{\alpha}^{\pm}k$.
\end{itemize}
\begin{proof}
If one of the two conditions hold, then the maximum in $p_{\alpha}^{\pm}(s)$ is always given by $\mp\mu(h_{\alpha})-d_\alpha sk$ for all $s\in\{i,i\pm1\}$ which implies equality in \eqref{tosim}. The general case follows from several case considerations and we illustrate the case  $\mp\mu(h_{\alpha})-d_\alpha (i+1)k\leq 0$ and $\mp\mu(h_{\alpha})-d_\alpha ik \geq 0$ only. In this case the inequality is equivalent to 
$$\mp2\mu(h_{\alpha})-2d_\alpha ik\leq \mp\mu(h_{\alpha})-d_\alpha (i-1)k \Leftrightarrow \mp\mu(h_{\alpha})\leq d_{\alpha}(i+1)k.$$
\end{proof}
\end{lem}
\subsection{} Now we prove several maps among the aforementioned modules and show in the subsequent section how they can be used to simplify their presentation. First we introduce some well-known sets which have been used already in \cite{CV13} and \cite{KV14} to simplify the relations of stable Demazure modules. Roughly speaking, our proof follows the proofs of \cite{CV13} and \cite{KV14} in a modified way. For $r\in\mathbb{N}$ and $s\in\bz_+$, let \begin{equation}\label{maxs}\bs(r,s)=\left\{(b_p)_{p\ge 0}: b_p\in\bz_+, \ \  \sum_{p\ge 0} b_p=r,\ \ \sum_{p\ge 0} pb_p=s\right\}. \end{equation}
For $k\in\bz_+$, we denote by $\bs(r,s)_k$ (resp. $_k\bs(r,s)$) the subset of $\bs(r,s)$ consisting of elements satisfying $$b_p=0,\ \ p\ge k,\ \ ({\rm{resp.}}\ \ b_p=0\ \ p<k).$$ 
Given $x\in\lie g$ define the following elements
\begin{align*}\mathbf{x}(r,s)_k&=\sum_{\mathbf{b}\in\bs(r,s)_k}(x\otimes 1)^{(b_0)}\cdots (x\otimes t^{k-1})^{(b_{k-1})},\\ _{k}\mathbf{x}(r,s)&=\sum_{\mathbf{b}\in _k\bs(r,s)}(x\otimes t^k)^{(b_k)}\cdots (x\otimes t^s)^{(b_s)}\\
~^t\!\mathbf{x}(r,s)&=
\sum_{\mathbf{b}\in\bs(r,s)}(x\otimes t)^{(b_0)}(x\otimes t^2)^{(b_1)}\cdots (x\otimes t^{s+1})^{(b_{s})},
\end{align*}
where for any integer $p$ and any $y\in\lie g[t]$,  we set $y^{(p)}=y^p/p!$. Moreover, we set for simplicity $_{0}\mathbf{x}(r,s)=\mathbf{x}(r,s)$ and define similarly as above the elements $~^t\!\mathbf{x}(r,s)_k$ and ${}^{t}_{k}\mathbf{x}(r,s)$.
We obtain,  \begin{equation}\label{special1} \mathbf{x}(1,s)=x\otimes t^s,\ \  ~^t\!\mathbf{x}(1,s)=x\otimes t^{s+1},\ \ _{k}\mathbf{x}(r,kr)=(x\otimes t^k)^{(r)}.\end{equation}
The following proposition is an easy modification of the analogues result proved in \cite{CV13}.
\begin{prop}\label{maxt} 
 Let $V$ be any representation of $\mathfrak{I}[t]$ and let $v\in V$, $x\in\lie g$, and $K,k\in\bz_+$. Then,
\begin{align*}{}^{t}\mathbf{x}(r,s)v=0\ {\rm{for\ all}}\ s\in\bz_+, r\in\bn
 \ \ {\rm{with}}\  s+r\ge 1+kr+K\ \iff\\  {}^{t}_{k}\mathbf{x}(r,s)v=0  \ {\rm{ for \ all}}\  s\in\bz_+,r\in\bn \  {\rm{with}}\ s+r\ge 1+ kr+K.\end{align*}
 Moreover, if $V$ is stable under the action with $(x\otimes 1)$, the same equivalence holds if we erase the $t$ in the superscript.
 \qed
\end{prop}
\subsection{} For a Lie algebra $\mathfrak{a}$ we denote by $\mathfrak{a}[t]_+$ the subspace spanned by the positive degree elements, i.e., $\mathfrak{a}[t]_+=\bigoplus_{r>0}\mathbb{C} (\mathfrak{a}\otimes t^r)$. Furthermore, let 
$$\mathbf{U}_{+}=\mathbf{U}(\lie g[t])\cdot \lie n^+[t]\oplus\mathbf{U}\big(\lie n^-[t]\oplus\lie h[t]\big)\cdot \lie h[t]_+$$
and define $\mathbf{U}_{++}$ similarly as above by replacing $\lie n^+[t]$ in $\mathbf{U}_+$ by $\lie n^+[t]_+$. Consider the automorphism $\Phi$ on $\lie g[t]$ given by 
$$x^+_\a\otimes t^s\mapsto x^-_\a\otimes t^s,\ \ x^-_\a\otimes t^s\mapsto x^+_\a\otimes t^s,\ \ s\geq 0$$ and note that
  $$\Phi(\mathbf{U}_{++})=\mathbf{U}(\lie g[t])\cdot \lie n^-[t]_+\oplus\mathbf{U}(\lie n^+[t]\oplus\lie h[t])\cdot \lie h[t]_+.$$
The following result is a slight modification of the Garland identities \cite[Lemma 7.5/7.8]{G78}. For a current algebra reformulation of the above lemma see also \cite[Lemma 1.3]{CP01}.
\begin{lem}\label{gar} {{Given  $s\in\bn $}}, $r\in\bz_+$ and $\alpha\in R^+$ we have the following identities
\begin{enumerate}
 \item $(x_{\alpha}^-\otimes t)^{(s)}(x_{\alpha}^+\otimes 1)^{(s+r)}\equiv (-1)^{s}\ \mathbf{x}^+_\a(r, s)
\ \ \mathrm{mod}\ \big(\Phi(\mathbf{U}_{++})\cap \mathbf{U}(\lie{I}[t])\big)$
 \vspace{0,2cm}
 
 \item $(x_{\alpha}^+\otimes 1)^{(s)}(x_{\alpha}^-\otimes t)^{(s+r)}\equiv (-1)^{s}\  {}^{t}\mathbf{x}^-_\a(r, s)
\ \ \mathrm{mod}\ \big(\mathbf{U}_+\cap \mathbf{U}(\lie{I}[t])\big)$
\end{enumerate}
 \qed
\end{lem}
\subsection{} The key result which will lead to the aforementioned simplified presentation is the following.
\begin{prop}\label{rel1} Let $\mu\in P$ and suppose that the tuple $_{\pm}\xi^{\alpha}$  from \eqref{part1} is a partition. Then, the following relations are true in the module $M'_{\mu,\mathbf{p}}$. For all $\alpha\in R^+(\mu)$ (resp. $\alpha\in R^-(\mu)$) and  $r\in\mathbb{N}, s, k\in \mathbb{Z}_+$ such that 
$$s+r\ge 1+ rk+p_{\alpha}^-(k+1) \ \ \ (\text{resp. } s+r\ge 1+ rk+p_{\alpha}^+(k))$$ we have ${}^{t}_{}\mathbf{x}^-_\alpha(r,s)v=0,\  (\text{resp. } \mathbf{x}^+_\alpha(r,s)v=0).$
\begin{proof} We give the proof of the proposition only for $\alpha\in R^+(\mu)$. The other case is treated similarly and we omit the details. For simplicity, we set in the rest of the proof $_{-}\xi^{\alpha}_j=\xi^{\alpha}_j$, $s_{\alpha}=s^-_{\alpha}$ and $p_{\alpha}=p^-_{\alpha}$. If $r\ge \xi_1^{\alpha}=\mu(h_{\alpha})-p_{\alpha}(1)$, then we get
$$s+r\ge 1+ k\xi_1^{\alpha}+p_{\alpha}(k+1)\ge 1+\mu(h_{\alpha})-\xi_{k+1}^{\alpha}\geq 1+p_{\alpha}(1)+\xi_{1}^{\alpha}-\xi_{k+1}^{\alpha}\ge p_{\alpha}(1)+1$$ and the claim follows from Lemma~\ref{gar} since $(x_{\alpha}^-\otimes t)^{a}v=0$ if $a\ge p_{\alpha}(1)+1$. So assume that $r<\xi^{\alpha}_1$. If additionally $r\le \xi^{\alpha}_{s_{\alpha}}=p_{\a}(s_{\alpha}-1)$ we can deduce $s+r\ge 1+(s_{\alpha}-1)r$ which we will explain now. To see the above inequality let $k<s_{\alpha}-1$ (otherwise the statement is clear) and we get
$$s+r\ge 1+kr+ p_{\alpha}(k+1)\ge 1+kr+(s_{\alpha}-k-1)\xi_{s_\alpha}^{\alpha}\ge 1+kr+(s_{\alpha}-k-1)r$$
and the claim follows. Moreover, we obtain that any $\mathbf{b}\in\bs(r,s)$ has the property $b_m>0$ for some $m\ge s_\alpha-1$, since otherwise we would have $$s=\sum_{p\le s_{\alpha}-2}pb_p\le r(s_{\alpha}-2).$$ In particular, the claim follows from the relation $(x_{\alpha}^-\otimes t^{s_{\alpha}})v=0$ in $M'_{\mu,\mathbf{p}}$. So we can assume in the rest of the proof that $\xi^{\alpha}_{s_{\alpha}}<r<\xi^{\alpha}_1$ and let $i\in\{1,\dots,s_{\alpha}-1\}$ be such that $\xi^{\alpha}_{i+1}<r\leq\xi^{\alpha}_i$.
Furthermore, let $\mathbf{b}\in\bs(r,s)_{s_{\alpha}-1}$ (otherwise we are done). Since $s=(s_{\alpha}-2)b_{s_{\alpha}-2}+\cdots +b_1$ we obtain 
 $$(s_{\alpha}-2)b_{s_{\alpha}-2}+\cdots +(i-1)b_{i-1}+(r-\sum_{j\geq i-1} b_{j})(i-2)\ge s\ge 1+ r(k-1)+p_{\alpha}(k+1),$$ and hence
$$(s_{\alpha}-i)b_{s_{\alpha}-2}+\cdots+b_{i-1}\ge 1+ r(k-i+1)+p_{\alpha}(k+1).$$
Since $r>\xi^{\alpha}_{i+1}$, we see that the right hand side of the inequality is bigger or equal to $p_{\alpha}(i)$ provided $k\ge i-1$:
\begin{align*}
1+ r(k-i+1)+&p_{\alpha}(k+1)\geq 1+(k-i+1)(p_{\alpha}(i)-p_{\alpha}(i+1))+p_{\alpha}(k+1)&\\&\geq 1+(p_{\alpha}(i)-p_{\alpha}(i+1))+\sum^{k-i}_{j=1} (p_{\alpha}(i+j)-p_{\alpha}(i+j+1))+p_{\alpha}(k+1)&\\&\geq 1+ p_{\alpha}(i).
\end{align*}
If $k< i-1$, then
$$(s_{\alpha}-i)b_{s_{\alpha}-2}+\cdots+b_{i-1}\ge 1+\sum_{i\geq  j\ge k+2}(\xi^{\alpha}_j-r)+\xi^\alpha_{i+1}+\cdots+\xi^\alpha_{s_{\alpha}}\geq 1+p_{\alpha}(i).$$
So each summand of ${}^{t}_{}\mathbf{x}^-_\alpha(r,s)$ acts as zero by the defining relations of $M_{\mu,\mathbf{p}}'$ and the calculations above. This finishes the proof of the proposition.
\end{proof}
\end{prop}
 \begin{cor}\label{maincor}Let $k\in\mathbb{N}$, $\alpha\in R^{\mp}(\mu)$ and $a_k,\dots, a_{s_{\alpha}^{\pm}}$ non-negative integers such that 
 $$a_k+2a_{k+1}+\cdots+(s_{\alpha}^{\pm}-k+1)a_{s_{\alpha}^{\pm}}\geq 1+ p_{\alpha}^{\pm}(k).$$
We set $$r=a_{k}+\cdots+a_{s_{\alpha}^{\pm}},\ \ s=(k-1)a_{k}+\cdots+(s_{\alpha}^{\pm}-1)a_{s_{\alpha}^{\pm}}.$$
We have the following relations in $M'_{\mu,\mathbf{p}}$:
$${}_{k-1}^{t}\mathbf{x}^-_\alpha(r,s)v=0 \ \ (\text{if $\alpha\in R^{+}(\mu)$}),\ \ \ {}_{k}^{}\mathbf{x}^{+}_\alpha(r,r+s)v=0 \ \ (\text{if $\alpha\in R^{-}(\mu)$}).$$ 
\begin{proof}The statement follows immediately from Proposition~\ref{maxt}  and Proposition~\ref{rel1}.
 \end{proof}
 \end{cor}
\begin{thm}\label{mmmr}
For $\mu\in P$ and $\mathbf{p}$ as in Proposition~\ref{rel1} we have 
$$M''_{\mu,\mathbf{p}}\twoheadrightarrow M'_{\mu,\mathbf{p}}\twoheadrightarrow   M_{\mu,\mathbf{p}}.$$
Moreover, the first map is an isomorphism if $_{\pm}\xi^{\alpha}_1=\cdots=_{\pm}\xi^{\alpha}_{s_{\alpha}^{\pm}-1}$ and the second map is an isomorphism if $_{\pm}\xi^{\alpha}_1\neq _{\pm}\xi^{\alpha}_{2}$.
\begin{proof}
The first surjective map follows immediately from Corollary~\ref{maincor} if we set $a_k=p_{\alpha}^-(k)+1$ and the remaining integers to be zero. Then, in the notation of the above Corollary, we have
$${}_{k-1}^{t}\mathbf{x}^-_\alpha(r,s)v=(x_{\alpha}^-\otimes t^k)^{p_{\alpha}^-(k)+1}v=0 \ \ (\text{if $\alpha\in R^{+}(\mu)$})$$
$$ {}_{k}^{}\mathbf{x}^{+}_\alpha(r,r+s)v=(x_{\alpha}^+\otimes t^k)^{p_{\alpha}^+(k)+1}v=0 \ \ (\text{if $\alpha\in R^{-}(\mu)$}).$$ 
The second surjectivity holds by definition. If $_{\pm}\xi^{\alpha}_1=\cdots=_{\pm}\xi^{\alpha}_{s_{\alpha}^{\pm}-1}$, it is clear that all relations in $M_{\mu,\mathbf{p}}'$ also hold in $M_{\mu,\mathbf{p}}''$. Now assume that $_{\pm}\xi^{\alpha}_1\neq _{\pm}\xi^{\alpha}_{2}$. To see that the second map is an isomorphism, it will be enough to show that the monomial relations from Definition~\ref{maindef}(3) (for fixed $i$) hold in $M_{\mu,\mathbf{p}}'$ when $ _{\pm}\xi_{i}^{\alpha}= \ _{\pm}\xi_{i+1}^{\alpha}$, i.e. $$2p_{\alpha}^{\pm}(i)=p_{\alpha}^{\pm}(i-1)+p_{\alpha}^{\pm}(i+1).$$
We will show this by downward induction on $i\in\{1,\dots,s_{\alpha}^{\pm}-1\}$. If $i=s_{\alpha}^{\pm}-1$, then the answer is immediate and the relation already holds in $M_{\mu,\mathbf{p}}''$. Otherwise, we can assume by induction that
$$a_i+a_{i+1}+\cdots+a_{s_{\alpha}^{\pm}}>p_{\alpha}^{\pm}(i)-p_{\alpha}^{\pm}(i+1)$$
and hence setting $a_{i-1}=0$ we get 
$$a_{i-1}+2a_i+3a_{i+1}+\cdots+(s_{\alpha}^{\pm}-i+1)a_{s_{\alpha}^{\pm}}\geq 1+2p_{\alpha}^{\pm}(i)-p_{\alpha}^{\pm}(i+1)=1+p_{\alpha}^{\pm}(i-1).$$
So either $ _{\pm}\xi_{i-1}^{\alpha}= \ _{\pm}\xi_{i}^{\alpha}$ or the relation holds in $M_{\mu,\mathbf{p}}'$. Continuing in this way we obtain the desired relation, since the process has to stop by the reason of $_{\pm}\xi^{\alpha}_1\neq _{\pm}\xi^{\alpha}_{2}$.
\end{proof}
\end{thm}
\begin{rem}
For general $\mu$ and $\mathbf{p}$ it is quite hard to figure out the connection between these three modules. If the aforementioned conditions do not hold, than these maps do not need to be isomorphisms. However, the surjectivity will be enough to simplify the presentation of all Demazure modules.  
\end{rem}
\section{Simplified presentations and applications: embeddings}\label{section4}
\subsection{}\label{section41}We first discuss the trivial consequence of Theorem~\ref{mmmr} and extend the results of \cite{CV13} from $\mathfrak{g}$-stable Demazure modules to the general setting. For $\mu\in P$ and $\alpha\in R^{\mp}(\mu)$ we denote by  $s^{\pm}_\alpha,m^{\pm}_\alpha\in\bz_+$ the unique integers so that  \begin{equation}\label{ggg1}\mp\mu(h_{\alpha})=(s^{\pm}_{\alpha}-1)d_{\alpha}k+m^{\pm}_{\alpha},\ \ \ 0<m^{\pm}_{\alpha}\le d_{\alpha}k\end{equation}
It follows that the function $\mathbf{p}$ from Example~\ref{mainex}(2) turns into
$$p_{\alpha}^{\pm}(i)=\mathrm{max}\{0, (s^{\pm}_{\alpha}-i-1)d_{\alpha}k+m^{\pm}_{\alpha}\},\ \ i\in\mathbb{Z}_+$$
In particular, from Lemma~\ref{n234} we obtain that $_{\pm}\xi^{\alpha}_i=_{\pm}\xi^{\alpha}_{i+1}$ for all $1\leq i\leq s_{\alpha}^{\pm}-2$ and  $_{\pm}\xi^{\alpha}_{s_{\alpha}^{\pm}-1}=_{\pm}\xi^{\alpha}_{s_{\alpha}^{\pm}}$ provided that  $m^{\pm}_{\alpha}=d_{\alpha}k$. So in this case, the module $M'_{\mu,\mathbf{p}}$ is generated by $v$ subject to the relations in Definition~\ref{maindef}(1)-(2) and for $\alpha\in R^{\mp}(\mu)$
$$(x_{\alpha}^{\pm}\otimes t^{s_{\alpha}^{\pm}})v=0$$
 $$(x_\alpha^{\pm}\otimes t^{s^{\pm}_\a-1})^{m^{\pm}_\a+1}v=0,\ \text{ if $m^{\pm}_{\alpha}<d_{\alpha}k$}$$
\begin{cor}\label{simdemee}
Let $\mu\in P$ and $k\in\mathbb{N}$. Then $\mathbf{D}_{\mu}^k$ is a cyclic $\mathbf{U}(\mathfrak{I}[t])$--module generated by a non-zero vector $v$ with the following relations: 
\begin{align*}
(h\otimes t^s)v=\delta_{s,0}\cdot \mu(h)v, \text{ for all } h\in \mathfrak{h}
\end{align*}
and for $\alpha\in R^{\mp}(\mu)$ we have
\begin{equation}\label{demazurerelations1}
   \big(x_\alpha^{\pm}\otimes t^{s^{\pm}_\a-1}\big)^{m^{\pm}_\a+1}v=0, \ \text{if $m^{\pm}_{\alpha}<d_{\alpha}k$}; \ \ \big(x_\alpha^{\pm}\otimes t^{s^{\pm}_\a}\big)v=0,
  \end{equation}
\begin{equation}\label{demazurerelations2}
 \left(x_\alpha^+\otimes \mathbb{C}[t]\right)v=0,\ \left(x_\alpha^-\otimes t\right)^{\mathrm{max}\{0,\ \mu(h_{\alpha})-d_\a  k\}+1}v=0,\  \text{ if $\alpha\in R^+(\mu)$}
  \end{equation}
\begin{equation}\label{demazurerelations3}
\left(x_\alpha^-\otimes t\mathbb{C}[t]\right)v=0,\  \left(x_\alpha^+\otimes 1\right)^{-\mu(h_{\alpha})+1}v=0,\  \text{ if $\alpha\in R^-(\mu)$}
 \end{equation}
\begin{proof}
Let us denote the module generated by $v$ with the simplified relations listed above by $\mathbf{N}_{\mu}^k$. The presentation of Demazure modules stated in Theorem~\ref{mathieu} gives a surjective map $\mathbf{N}_{\mu}^k\twoheadrightarrow \mathbf{D}_{\mu}^k$. In order to prove that the map is an isomorphism, we only have to show that the defining relations of $\mathbf{D}_{\mu}^k$ hold in $\mathbf{N}_{\mu}^k$. But $\mathbf{N}_{\mu}^k$ is isomorphic to a quotient of $M'_{\mu,\mathbf{p}}\cong M''_{\mu,\mathbf{p}}$ and hence all relations also hold in $\mathbf{D}_{\mu}^k$. 
\end{proof}
\end{cor}
\begin{rem} In the special case $k=1$, a similar observation as in \cite[Proposition 3.4]{CV13} shows that some of the relations in \eqref{demazurerelations1} are still redundant. To be more precise, the second relation in \eqref{demazurerelations1} is only needed if $d_{\alpha}>1$ and the first relation in \eqref{demazurerelations1} is only needed when $d_{\alpha}=3=m^{\mp}_{\alpha}+2$. For example, this is one of the key observations why $\mathbf{D}_{\mu}^1$ for anti-dominant weights can be identified with local Weyl modules in the simply-laced case. 
\end{rem}
So given a module $L$ by generators and relations, let $L^{\mathrm{simp}}$ a module obtained by erasing certain relations in the presentation. In particular we have $L^{\mathrm{simp}}\rightarrow L \rightarrow 0$. The key idea above for showing that $L^{\mathrm{simp}}\cong L$ was to identify $\mu$ and $\mathbf{p}$ (the induced tuple in \eqref{part1} has to be a partition) such that there exists surjections
$$M_{\mu,\mathbf{p}}'\rightarrow L^{\mathrm{simp}}\rightarrow 0 \ \text{ and } \ L\rightarrow M_{\mu,\mathbf{p}}''\rightarrow 0.$$
This strategy can be used to simplify the presentation of certain more modules $V(\xi)$ introduced in \cite{CV13}. On can read the isomorphism $L^{\mathrm{simp}}\cong L$ also conversely, namely that we have more relations in $L^{\mathrm{simp}}$. This can be used for example to show that certain relations hold in the generalized Weyl module \cite{FM17a}.
\subsection{} The simplified presentation of Demazure modules allow us to embed $\mathbf{D}_{\mu}^{rk}$ into a $k$-fold tensor product of Demazure modules of the form $\mathbf{D}_{\mu'}^{r}$ (see Theorem~\ref{embed1}). In the limit $k\rightarrow \infty$ this will become for anti-dominant weights the well-known embedding of irreducible $\mathfrak{g}$-modules into a tensor product of fundamental irreducible $\lie g$-modules. The interest of this embedding is twofold. 

On the one hand we hope that this embedding can be used to come up with explicit combinatorial models for the Demazure crystal of $\mathbf{D}_{\mu}^k$ viewing it as a ``connected component'' inside the tensor product of suitable lower level crystals.

On the other hand, the motivation comes from quantum affine algebras and graded limits of irreducible finite-dimensional representations. It is sometimes important to prove surjective maps from the graded limit of a suitable representation to a higher level Demazure module as pointed out, for example, for level two in the article  \cite{BCMo15} for the graded limit of a prime irreducible object in the Hernandez-Leclerc category. To achieve such a map, the authors proved in type $A$ and $k=2$ (see \cite[Theorem 4]{BCMo15}) an embedding of Demazure modules. The aim of this section is to generalize their result.\\\

We state the main theorem of this section. For $k\in\mathbb{N}$ and $\mu\in P$ we define  
$$\mathcal{P}(\mu,k)=\{(\mu_1,\dots,\mu_{k})\in P^{\times {k}}: \mu_1+\cdots+\mu_{k}=\mu\}.$$
We say that an element $\boldsymbol{\mu}=(\mu_1,\dots,\mu_{k})\in \mathcal{P}(\mu,k)$ is \textit{pre-admissible} if
\begin{equation}\label{preadmissible}\alpha\in R^{\pm}(\mu)\Rightarrow \alpha\in R^{\pm}(\mu_i),\ \  1\leq i\leq k\end{equation}
and for such an element $\boldsymbol{\mu}$ we define non-negative integers $t_{\alpha}$ and $p^{\alpha}_j(\boldsymbol{\mu}),\ 0\leq j\leq t_{\alpha}$, $\alpha\in R^{\mp}(\mu)$ by the following equation
\begin{equation}\label{reord}\big(\mp\mu_{i_1}(h_{\alpha}),\dots, \mp\mu_{i_k}(h_{\alpha})\big)=\big(x,\dots,x,(x-1),\dots,(x-1),\dots,(x-t_{\a}),\dots,(x-t_{\a})\big)\end{equation}
where $p_j^{\alpha}(\boldsymbol{\mu})$ encodes how often $(x-j)$ appears in the above tuple. By convention we assume that we have $p_0^{\alpha}(\boldsymbol{\mu})>0$. Given $r\in \mathbb{N}$ we define also an element $0<m^{\pm}_{\alpha}(\boldsymbol{\mu},r)\leq d_{\alpha}r$ determined by 
$$x=\mp\mu_{i_1}(h_{\alpha})= m^{\pm}_{\alpha}(\boldsymbol{\mu},r) \mod d_{\alpha}r$$
We say that a pre-admissible element $\boldsymbol{\mu}$ is \textit{r-admissible} if it satisfies 
\begin{equation}\label{admissible}
m^{\pm}_{\alpha}(\boldsymbol{\mu},r)\cdot k> \sum_{j=0}^{t_{\a}} j\cdot  p^{\alpha}_j(\boldsymbol{\mu})\ \ \forall \alpha\in R^{\mp}(\mu)
\end{equation}
and
\begin{equation}\label{extra}
\mu(h_{\alpha})>k d_{\alpha}r \Rightarrow x=\mu_{i_1}(h_{\alpha})\geq (t_{\a}+d_{\alpha}r),\ \ \forall \alpha\in R^+(\mu)
\end{equation}
We emphasize that \eqref{extra} is only needed for roots in $R^+(\mu)$. It is not hard to see that each tuple $(\mu_1,\dots,\mu_k)$ is $r$-admissible for some $r>>0$. So it is natural to ask for the minimal $r$ such that a tuple is $r$-admissible. We will discuss this question in the next subsection and conjecture that there exists always a 1-admissible element. The main theorem of this section is the following.
\begin{thm}\label{embed1}
Given $\mu\in P$ and $r,k\in \mathbb{N}$, let  $\boldsymbol{\mu}=(\mu_1,\dots,\mu_k)\in\mathcal{P}(\mu,k)$ be $r$-admissible. Then  we have an injective map of graded $\mathfrak{I}[t]$-modules
$$\mathbf{D}_{\mu}^{rk}\hookrightarrow \mathbf{D}_{\mu_1}^{r}\otimes \cdots \otimes \mathbf{D}_{\mu_k}^{r},\ \ v_{\mu}\mapsto v_{\mu_1}\otimes \cdots \otimes v_{\mu_k}.$$
\end{thm} 
The proof of this theorem occupies the rest of this section.
\subsection{}The fact that $\mathbf{D}_{\mu}^{k}$ appears as a subquotient 
follows from the next result which can be derived from \cite{Ku88a}.
\begin{prop} Let $r\in\mathbb{N}$ and choose positive integers $k_1,\dots,k_r$ such that $k_1+\cdots+k_r=k$ and $(\mu_1,\dots,\mu_r)\in \mathcal{P}(\mu,r)$. Then we have a surjective map 
$$\mathbf{D}_{\mu_1}^{k_1}\otimes \cdots \otimes \mathbf{D}_{\mu_k}^{k_r}\supseteq \mathbf{U}(\lie I[t])(v_{\mu_1}\otimes \cdots \otimes v_{\mu_r})\rightarrow \mathbf{D}_{\mu}^k\rightarrow 0$$
\hfill\qed
\end{prop}
The previous proposition implies that it is enough to show that for an $r$-admissible element $(\mu_1,\dots,\mu_k)\in \mathcal{P}(\mu,k)$ the weight vector $(v_{\mu_1}\otimes\cdots \otimes v_{\mu_k})$ satisfies the defining relations of $\mathbf{D}_{\mu}^{rk}$ which were simplified in Corollary~\ref{simdemee}. Without this simplification it was not possible for us to show that the required relations hold. The following lemma is trivial but useful:
\begin{lem}\label{techcon1}
Let $x,j\in \mathbb{Z}_+$ such that $x\geq j$ and $\ell\in\mathbb{N}$. We write $x=(s-1)\ell+m$,\  $x-j=(s'-1)\ell+m'$ and $j=(q-1)\ell+m''$, $0<m,m',m''\leq \ell$. Then we have
 $$(s, m)=
 \begin{cases}
 (s'+q-1, m'+m'') & \ \text{if $m'+m''\le \ell$} \\
 (s'+q, m'+m''-\ell) & \ \text{if $m'+m''> \ell$}\\
\end{cases}$$
\hfill\qed
\end{lem}
\subsection{Proof of Theorem~\ref{embed1}} Fix an element $\alpha\in R^{\mp}(\mu)$. As usual we define the numbers $s_{\alpha}^{\pm}, m_{\alpha}^{\pm}$ for $\mu$ as in \eqref{ggg1} (with respect to $d_{\alpha}kr$) and denote the corresponding elements for $\mu_i$ by  $(s^i_{\alpha})^{\pm}, (m^i_{\alpha})^{\pm}$ (with respect to $d_{\alpha} r$). Recall that $\boldsymbol{\mu}$ satisfies property \eqref{preadmissible}. In the rest of this section we abbreviate $v_j$ for the cyclic generator of $\boldsymbol{D}_{\mu_j}^r$, $1\le j\le k$ and let $v$ be the cyclic generator of $\boldsymbol{D}_{\mu}^{rk}$. We abbreviate further $p_j=p_j^{\alpha}(\boldsymbol{\mu})$ and $t=t_{\alpha}$ and recall that $p_0+\cdots +p_t=k.$ We have
\begin{align}\notag(s_{\alpha}^{\pm}-1) d_\alpha kr+m_{\alpha}^{\pm}&=\mp \mu(h_\alpha)=p_0x+p_1(x-1)+\cdots +p_t(x-t)&\\&\notag=kx-(p_1+2p_2+\cdots+t p_t)&\\&\label{hhhtt}
=\big((s^{i_{1}}_{\alpha})^{\pm}-1\big)d_\alpha kr+\big(km^{\pm}_{\alpha}(\boldsymbol{\mu},r)-(p_1+2p_2+\cdots+t p_t)\big).\end{align}
From \eqref{admissible} and $0<m^{\pm}_{\alpha}(\boldsymbol{\mu},r)\le d_\alpha r$ we know that 
$$0<km^{\pm}_{\alpha}(\boldsymbol{\mu},r)-(p_1+2p_2+\cdots+t p_t)\leq d_\alpha r$$ and hence $s_{\alpha}^{\pm}=(s^{i_{1}}_{\alpha})^{\pm}$. Now using Lemma~\ref{techcon1} we also get that $$s_{\alpha}^{\pm}=(s^{i_{1}}_{\alpha})^{\pm}\ge (s^{j}_{\alpha})^{\pm},\text{ for all $1\le j\le k$}$$
and hence we have proved the connection
 \begin{equation}\label{conneq1}s_{\alpha}^{\pm}=\max\{ (s^i_{\alpha})^{\pm}: 1\leq i\leq k\}.\end{equation}
\textit{Relations in \eqref{demazurerelations1}}: The relation $(x_\alpha^{\pm}\otimes t^{s^{\pm}_\a})(v_1\otimes \cdots \otimes v_k)=0$ is immediate from \eqref{conneq1} and Corollary~\ref{simdemee}. So we argue in the remaining part why the following relation holds:
   $$(x_\alpha^{\pm}\otimes t^{s^{\pm}_\a-1})^{m_\alpha^{\pm}+1}(v_1\otimes \cdots \otimes v_k)=0.$$
So assume by contradiction that the above relation does not hold. Write $Y=(x_\alpha^{\pm}\otimes t^{s^{\pm}_\a-1})$ and $m=m_\alpha^{\pm}+1$, then
$$Y^{m}(v_1\otimes \cdots \otimes v_k)=
\sum\limits_{j_1+\cdots +j_k=m}Y^{j_1}v_1\otimes \cdots \otimes Y^{j_k}v_k.$$
Note that we either have $Yv_i=0$ (which is the case if $s^{\pm}_\a\ge (s^{i}_{\alpha})^{\pm}+1$) or $s^{\pm}_\a\leq (s^{i}_{\alpha})^{\pm}$ and hence by \eqref{conneq1} we must have equality $s^{\pm}_\a=(s^{i}_{\alpha})^{\pm}$. In the latter case we have
$Y^{s}v_i=0$ if $s\ge (m^{i}_{\alpha})^{\pm}+1$. So there must exist a tuple $(j_1,\dots, j_k)$ in the above sum which has the following conditions:
\begin{enumerate}
    \item  $j_1+\cdots +j_k=m$
    \item $j_i=0$ if  $s^{\pm}_\a \neq (s^{i}_{\alpha})^{\pm}$
    \item $j_i\le (m^{i}_{\alpha})^{\pm}$ if  $s^{\pm}_\a = (s^{i}_{\alpha})^{\pm}$
\end{enumerate}
Let $A=\{1\leq i\leq k : s^{\pm}_\a = (s^{i}_{\alpha})^{\pm}\}$ and set $I_j=\big[p_0+\cdots+p_{j-1}+1,p_0+\cdots+p_{j}\big]$. From Lemma~\ref{techcon1} we see that an element $i\in I_j$ lies in $A$ if and only if 
$$0\leq j\leq d_{\alpha} r\ \ \text{and}\ \  (m_{\alpha}^i)^{\pm}+j\leq d_{\alpha}r.$$
In particular, if $i\in I_j\cap A$ we have the equality $(m^{i_{1}}_{\alpha})^{\pm}=(m^{i}_{\alpha})^{\pm}+j$ and therefore
\begin{equation}\label{zzttrr}m=\sum\limits_{i\in A}j_i=\sum_{j=0}^t\ \sum_{i\in A\cap I_j}j_i\le \sum_{j=0}^t |A\cap I_j|\big((m^{i_{1}}_{\alpha})^{\pm}-j\big)\end{equation}
On the other hand from \eqref{hhhtt} we get
$$m=k(m^{i_{1}}_{\alpha})^{\pm}+1-(p_1+2p_2+\cdots+t p_t).$$
So if we substitute this into inequality \eqref{zzttrr} we end up in a contradiction since $|A\cap I_j|\leq p_j$:
$$p_0(m^{i_{1}}_{\alpha})^{\pm}+p_1\big((m^{i_{1}}_{\alpha})^{\pm}-1\big)+\cdots +
p_t\big((m^{i_{1}}_{\alpha})^{\pm}-t\big)< \sum_{j=0}^t |A\cap I_j|\big((m^{i_{1}}_{\alpha})^{\pm}-j\big).$$

\textit{Relations in \eqref{demazurerelations2}}: Now suppose that $\alpha\in R^{+}(\mu)$. The only non-trivial relation is given by
\begin{equation}\label{fffdd}(x_\alpha^-\otimes t)^{\mathrm{max}\{ 0, \mu(h_\a)-d_\a kr \}+1}(v_1\otimes \cdots \otimes v_k)=0. \end{equation}
Suppose that $\mu(h_\a)\le kd_\a r$ then $s_\alpha^-\leq 1$ which implies together with \eqref{conneq1} that $(s^{i}_{\alpha})^-\leq 1$ for each $1\le i\le k$. Hence the above relation is immediate in this case. Suppose from now on that $\mu(h_\a) > k d_\a r$ and note that $\mu_i(h_{\alpha})\geq d_{\alpha}r$ for all $1\leq i \leq k$ by \eqref{extra}. Again we write $Y=(x_\alpha^-\otimes t)$ and consider a typical summand on the left hand side of \eqref{fffdd}
$$Y^{j_1}v_1\otimes \dots \otimes Y^{j_k}v_k.$$
So if for all $i\in\{1,\dots,n\}$ the inequality $j_i\leq \max\{0,\mu_i(h_{\alpha})-d_{\alpha}r\}$ holds we would get a contradiction togehter with \eqref{extra}:
$$\mu(h_\a) - k d_\a r+1=j_1+\dots+j_k\leq \sum_{i=1}^k\max\{0,\mu_i(h_{\alpha})-d_{\alpha}r\}\leq  \mu(h_\a) - k d_\a r$$
\textit{Relations in \eqref{demazurerelations3}}: The first one follows directly from property \eqref{preadmissible} and the second relation is standard. 
This completes the proof.

\subsection{} Given $\mu\in P$ and $k\in\mathbb{N}$, there are a few natural candidates for $\boldsymbol{\mu}$ being $1$-admissible arising from \cite{F14aa} which we will explain now. Let $\sigma\in W$ such that $\sigma(\mu)\in P^+$ and consider
$$\mathcal{P}^+(\sigma(\mu),k)=\{(\lambda_1,\dots,\lambda_{k})\in (P^+)^{\times {k}}:\lambda_1+\cdots+\lambda_{k}=\sigma(\mu)\}.$$
This set is partially ordered and has a unique maximal element $\boldsymbol{\lambda}^{\mathbf{max}}$ maximizing the dimension of tensor products of the corresponding irreducible representations \cite{F14aa}. This unique element can be explicitly defined (see \cite[Section 3.1]{F14aa}) and gives at least a pre-admissible tuple 
\begin{equation}\label{fouasd}\boldsymbol{\mu}=(\sigma^{-1}(\lambda_1),\dots,\sigma^{-1}(\lambda_k))\in \mathcal{P}(\mu,k),\ \ \boldsymbol{\lambda}^{\mathbf{max}}=(\lambda_1,\dots,\lambda_k).\end{equation}
However this element is generically only $1$-admissible in type A which we illustrate in the next example.
\begin{example}\label{exfou}
\begin{enumerate}
\item Let $\lie g$ be of type A and $\mu\in P$, $\sigma\in W$ and $\boldsymbol{\lambda}^{\mathbf{max}}\in \mathcal{P}^+(\sigma(\mu),k)$ as above. From the explicit description (see \cite[Section 3.1]{F14aa}) we see that $\boldsymbol{\lambda}^{\mathbf{max}}$ satisfies  
\begin{equation}\label{kl231}|\lambda_i(h_{\alpha})-\lambda_j(h_{\alpha})|\leq 1,\ \ \forall \alpha\in R^+.\end{equation}
Hence the element $\boldsymbol{\mu}=(\sigma^{-1}(\mu_1),\dots,\sigma^{-1}(\mu_k))$ is $r$-admissible for all $r\geq 1$ (we have $t_{\alpha}\leq 1$ in this case). 
\item Let $\mathfrak{g}$ be of type $C_2$, $k=2$, and $\mu=2\varpi_1+\varpi_2$. The unique maximal element from \cite{F14aa} is given by $(\varpi_1+\varpi_2,\varpi_1)$, however \eqref{admissible} is violated for $r=1$ and $\alpha=\alpha_1+\alpha_2$. So it is not 1-admissible, but 2-admissible. Nevertheless, $(2\varpi_1,\varpi_2)$ provides a 1-admissible element which is also 2-admissible. 
\end{enumerate}
\end{example}
So the previous example shows that in type A the unique maximal element from \cite{F14aa} gives rise to a 1-admissible element in $\mathcal{P}(\mu,k)$.  We conjecture that there is always a $1$-admissible element. For now, we can only prove the following.
\begin{lem}
Let $\boldsymbol{\mu}$ be the pre-admissible element from \eqref{fouasd} and let $\lie g$ be of non-exceptional type. Then $\boldsymbol{\mu}$ is $1$-admissible or there exists $\alpha\in R^{\mp}(\mu)$ such that
$t_{\alpha}=2=(m_{\alpha}^{i_1})^{\pm}+1$. If additionally $$\mp\mu_{i_1}(h_{\alpha})\equiv (d_{\alpha}+1) \mod 2d_{\alpha}$$
for those $\alpha$, then $\boldsymbol{\mu}$ is $2$-admissible.
\begin{proof}
For $\lie g$ of non-exceptional type we always have $t_{\alpha}\leq 2$ since \eqref{kl231} holds for all roots $\alpha\in R^+$ such that $\varpi_j(h_{\alpha})\leq 1$ for all $j\in I$ (see \cite[Section 3.1]{F14aa}) and every coroot can be written as a sum of two such coroots. Assume that $\boldsymbol{\mu}$ is not  $1$-admissible. Since $t_{\alpha}\leq 1$ for all $\alpha\in R^+$ would imply $1$-admissibility, we deduce that there exists a root $\alpha\in R^+$ satisfying $t_{\alpha}=2$. 
Among those $\alpha$ with $t_{\alpha}=2$ there must exist at least one which satisfies $(m_{\alpha}^{i_1})^{\pm}= 1$; otherwise we would again get $1$-admissibility. This show the first part of the statement. Now let $\alpha\in R^{\mp}(\mu)$ such that $t_{\alpha}=2=(m_{\alpha}^{i_1})^{\pm}+1$ and write
$$\mp \mu_{i_1}(h_{\alpha})=q_{\alpha} d_{\alpha}+(m_{\alpha}^{i_1})^{\pm}=q'_{\alpha} 2d_{\alpha}+(d_{\alpha}+1)$$
where the second equality follows from the additional assumption. So the new remainder is $(d_{\alpha}+1)$ with respect to $2d_{\alpha}$ and hence greater or equal to $2$.
Now if we have an $\alpha\in R^{\mp}(\mu)$ with $t_{\alpha}= 2$ and $(m_{\alpha}^{i_1})^{\pm}\geq 2$, the new remainder can be calculated as follows. We write
$$\mp \mu_{i_1}(h_{\alpha})=q_{\alpha} d_{\alpha}+(m_{\alpha}^{i_1})^{\pm}$$
as above. If $q_{\alpha}\equiv 0 \mod 2$, the remainder with respect to $2d_{\alpha}$ does not change or $q_{\alpha}\equiv 1 \mod 2$ and the new remainder is 
$(m_{\alpha}^{i_1})^{\pm}+d_{\alpha}\leq 2 d_{\alpha}.$ In both case we increase the remainder and we have proved the following. For all $\alpha\in R^{\mp}(\mu)$ we have $t_{\alpha}\leq 1$ or $t_{\alpha}=2$ in which case the remainder with respect to $2d_{\alpha}$ is greater or equal to $2$. This shows that \eqref{admissible} and \eqref{extra} hold for $r=2$ and $\boldsymbol{\mu}$ is $2$-admissible.
\end{proof}
\end{lem}

\subsection{} We assume in this subsection that $\lie g$ is of simply-laced type. As mentioned earlier, we can give a lower bound for the dimension of certain representations of the quantum affine algebra $\widehat{\mathbf{U}}_q(\lie g)$. We denote by $\mathcal{P}^+$ the monoid generated by $\boldsymbol{\varpi}_{i,a}$ where $i\in I$ and $a\in\mathbb{C}(q)^{\times}$ and let $\mathcal{P}_{\mathbb{Z}}^+$ the submonoid generated by the elements $\boldsymbol{\varpi}_{i,a}$ with $a\in q^{\mathbb{Z}}$. It is well-known \cite{CP91a,CP95} that the isomorphism classes of irreducible finite-dimensional representations of $\widehat{\mathbf{U}}_q(\lie g)$ are indexed by elements in $\mathcal{P}^+$. For $\boldsymbol{\pi}\in \mathcal{P}_{\mathbb{Z}}^+$ we denote the corresponding finite-dimensional irreducible module by $V(\boldsymbol{\pi})$ and its graded limit by $L(\boldsymbol{\pi})$ (for a precise definition and comments on the existence see \cite[Section 1.8]{BCMo15} for example). Note that we have a weight map $\mathrm{wt}:\mathcal{P}^+\rightarrow P^+$ given by extending the assignment $\mathrm{wt}(\boldsymbol{\varpi}_{i,a})=\varpi_i$ to a morphism of monoids.
\begin{cor} Let $\lie g$ be of simply-laced type. 
Assume that $\boldsymbol{\pi}^1,\dots,\boldsymbol{\pi}^{k}\in \mathcal{P}_{\mathbb{Z}}^+$ are such that
\begin{enumerate}
    \item There exists a map of $\widehat{\mathbf{U}}_q(\lie g)$-modules $$V(\boldsymbol{\pi})\rightarrow V(\boldsymbol{\pi}^1)\otimes \cdots \otimes V(\boldsymbol{\pi}^k),\ \ \boldsymbol{\pi}=\boldsymbol{\pi}^1\cdots \boldsymbol{\pi}^k$$
     \item $V(\boldsymbol{\pi}^i)$ is a quantum Weyl module for all $1\leq i\leq k$ \vspace{0,1cm}
    \item The element $(w_0\mathrm{wt}(\boldsymbol{\pi}^1),\dots,w_0\mathrm{wt}(\boldsymbol{\pi}^k))$ is $1$-admissible.
\end{enumerate}
Then $\dim V(\boldsymbol{\pi})\geq \dim \mathbf{D}_{w_0\mathrm{wt}(\boldsymbol{\pi})}^k$.
\begin{proof}
From property (1) and \cite[Lemma 2.20]{Mou10} we know that there exists a map $L(\boldsymbol{\pi})\rightarrow L(\boldsymbol{\pi}^1)\otimes \cdots \otimes L(\boldsymbol{\pi}^k)$ on the level of graded limits which is surjective onto the image. By property (2) we obtain an isomorphism between $L(\boldsymbol{\pi}^i)$ and $\mathbf{D}_{w_0\mathrm{wt}(\boldsymbol{\pi}^i)}^1$. Hence the image is isomorphic to $\mathbf{D}_{w_0\mathrm{wt}(\boldsymbol{\pi})}^k$ by property (3) and Theorem~\ref{embed1}. The rest follows from the fact that the graded limit preserves dimensions.
\end{proof}
\end{cor}
Although the conditions of the above corollary might look strong, there are sufficient conditions worked out for property (1) and (2) to hold (see for example \cite{Ch002,CP01}). Using the formula for the denominator of the normalized $R$-matrix (see \cite{Fu22}) and Hernandez reduction theorem \cite{He09} there are even necessary and sufficient conditions for property (2). 
\subsection{}

It is still an open problem to determine the classical decompositions of Demazure modules with respect to a maximal semi-simple Lie subalgebra $\lie g_0\subseteq \lie g$. The problem in the $\lie g$-stable case has some solutions and is easier to understand for rectangular anti-dominant weights (i.e. $\mu\in -k\cdot P^{+}$); see for example \cite{C01,FoL06}. However, the question seems difficult for prime Demazure modules (see \cite{BK20b,Sh15a} for some progress). In the $\lie g$-stable case there is another way to attack this problem using tensor products of Kirillov-Reshetikhin crystals. 

If $B^{r,s}$ denotes the KR crystal of classical highest weight $s\varpi_r$ and $B=\bigotimes_{j=1}^{N} B^{r_j,s_j}$ is a tensor product of KR crystals bounded by $k$, i.e. $\lceil \frac{s_j}{c_{r_j}} \rceil \leq k$ for all $1\leq j \leq N$ (see \cite{ST12} for the table defining the integers $c_r$), then the proof of \cite[Theorem 3.8]{LS19a} shows for nonexceptional types that the Demazure crystal of $\mathbf{D}_{\mu}^k$ for anti-dominant $\mu=\sum_{j=1}^N r_jw_0\varpi_j$ can be obtained by the following steps. 

$\bullet$ Remove all edges in $B$ which are contained in the length $k$ head of a $0$-string; following \cite{LS19a} we call this subcrystal by $\tilde{D}_k(B)$.\vspace{0,1cm}

$\bullet$ Then, up to certain $0$-arrows, the Demazure crystal of $\mathbf{D}_{\mu}^k$ is the connected component $D_k(B)$ of $\tilde{D}_k(B)$ containing the unique element of weight $\mu$.

Removing all $0$-arrows of $D_k(B)$ in a third step, provides the $\lie g$-decomposition of $\mathbf{D}_{\mu}^k$. Finding very explicit combinatorial models for these subcrystals is another problem. 

Now if $\lie g$ and $\mu$ are arbitrary such that $R^-(\mu)\neq \emptyset$, we consider the maximal semi-simple Lie subalgebra $\mathfrak{g}_0\subseteq \lie g$ such that the module $\mathbf{D}_{\mu}^k$ is $\lie g_0$-stable. An interesting problem is to describe explicitly the Demazure crystal of $\mathbf{D}_{\mu}^k$ as the connected component in the corresponding tensor product of Demazure crystals containing the unique element of weight $\mu$. We emphasize here that the way of splitting $\mu$ is very important. We will address this problem and precise statements elsewhere and discuss only an example. 
\begin{example}
Let $\lie g$ be of type $A_2$ and $\mu=\varpi_1-2\varpi_2$. Then $$V_{s_2s_1}(\Lambda_1+\Lambda_2)=\mathbf{D}^2_{\mu}\hookrightarrow \mathbf{D}^1_{-\varpi_2}\otimes \mathbf{D}^1_{\varpi_1-\varpi_2}=V_{s_2s_1}(\Lambda_1)\otimes V_{s_2}(\Lambda_2).$$
We know from \cite{Ka93a} that the crystal base $B_w(\Lambda)$ of $V_w(\Lambda)$ is the full subgraph of the crystal base of $V(\Lambda)$ whose vertices consist precisely of those elements that are reachable by raising operators from the unique element of weight $w\Lambda$. In this particular case we have, as sets,
$$B_{s_2s_1}(\Lambda_1)=\big\{v_{\Lambda_1}, \tilde{f}_1v_{\Lambda_1},\tilde{f}_2\tilde{f}_1v_{\Lambda_1}\big\}\ \ B_{s_2}(\Lambda_2)=\{v_{\Lambda_2},\tilde{f}_2 v_{\Lambda_2}\big\}$$
Now, tensoring the corresponding crystals 
$$
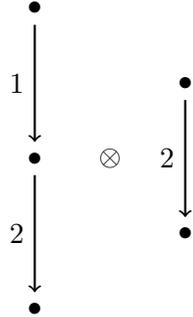
\begin{tikzpicture}
    \node (A) at (0,4) {$\bullet$};
    \node (B) at (0,2) {$\bullet$};
     \node (C) at (0,0) {$\bullet$};
          \node (D) at (1,2) {$\otimes$};
          \node(E) at (2,3) {$\bullet$};
           \node(F) at (2,1) {$\bullet$};
    \draw[thick,->] (A) -- (B) node[midway,left,rotate=0] {$1$};
       \draw[thick,->] (B) -- (C) node[midway,left,rotate=0] {$2$};
        \draw[thick,->] (E) -- (F) node[midway,left,rotate=0] {$2$};
  \end{tikzpicture}
$$
and considering the connected component containing the element of weight $\mu$ gives the graph
$$
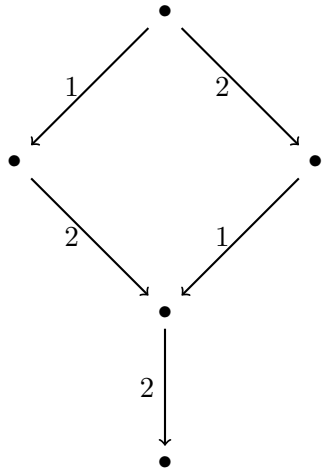
\begin{tikzpicture}
    \node (A) at (0,4) {$\bullet$};
     \node (B) at (-2,2) {$\bullet$};
      \node (C) at (2,2) {$\bullet$};
       \node (E) at (0,0) {$\bullet$};
        \node (F) at (0,-2) {$\bullet$};
         \draw[thick,->] (A) -- (B) node[midway,left,rotate=0] {$1$};
          \draw[thick,->] (A) -- (C) node[midway,left,rotate=0] {$2$};
           \draw[thick,->] (B) -- (E) node[midway,left,rotate=0] {$2$};
            \draw[thick,->] (C) -- (E) node[midway,left,rotate=0] {$1$};
             \draw[thick,->] (E) -- (F) node[midway,left,rotate=0] {$2$};
    \end{tikzpicture}
$$
Deleting all arrows whose label is not $2$ gives the $\mathfrak{sl}_2$ decomposition (corresponding to the simple root $\alpha_2$) of $\mathbf{D}^2_{\mu}$. In particular, $\dim (\mathbf{D}^2_{\mu})=5$. 
\end{example}

\bibliographystyle{plain}
\bibliography{bibfile}

\end{document}